%% file: arXivhomo.tex
\documentclass[11pt]{article}
\usepackage {epsfig} 
\usepackage{amsmath} 
\usepackage{amssymb}
\usepackage{latexsym} 
\usepackage{theorem}
\input{preamble4}

\def\mylabel#1{\label{#1}}
\reversemarginpar   
\begin{document}

 \title {\bf Periodic Orbits and Homoclinic Loops for Surface
Homeomorphisms}
\author
{
\Large  Morris W. Hirsch\thanks
{This research was supported in part by the
National Science Foundation.}\\
University of California at Berkeley\\[1ex]
University of Wisconsin at Madison
}
\maketitle

{\small This is a corrected version of the article in
  Michigan Mathematical Journal {\bf 47} (2000), no. 2, 395--406.}

\begin{abstract}
Let $p$ be a saddle fixed point for an orientation-preserving surface
diffeomorphism $f$, admitting a homoclinic point $p'$.  Let $V$ be an
open 2-cell bounded by a simple loop formed by two arcs joining $p$ to
$p'$, lying respectively in the stable and unstable curves at $p$.
It is shown that $f|V$ 
has fixed point index $\rho \in \{1,2\}$
where
$\rho$ depends only on the geometry of $V$ near $p$.
More generally,
for every $n\ge 1$, the union of the $n$-periodic orbits in $V$ is a
block of fixed points for $f^n$ whose index is $\rho$.

\end{abstract}

\maketitle

 \section{Introduction}\mylabel{sec:intro}


Poin\-car\'e invented homoclinic orbits, conjectured their existence
in the planar three body problem, and despaired of understanding their
complexity.  Research by Birk\-hoff, Cartwright \& Littlewood, and
Levinson revealed that near transverse homoclinic points there are
robust periodic points.  A pinnacle of this line of research, and the
basis for much of modern dynamical theory, is Smale's ``horseshoe''
theorem \cite{Sm65}.  For a diffeomorphism $f$ of a manifold of any
dimension, it states that every neighborhood of a transverse
homoclinic point meets a structurally stable, hyperbolic compact
invariant set $K$ on which some iterate $f^k$ is topologically
conjugate to the shift map on the Cantor set $2^{\ZZ}$.

Similar results have been obtained under weakenings of the
transversality assumption, including work by Burns \& Weiss
\cite{BW95}, Collins \cite{Co99}, Churchill \& Rod \cite{CR80},
Gavrilov \& \v{S}ilnikov \cite{GS72, GS73}, Guckenheimer \& Holmes
\cite{GH83}, Mischaikow \cite{Mi95}, Mischaikow \& Mrozek \cite{MM95},
Newhouse \cite{Ne80}, Rayskin \cite{Rayskin05}.

Among many important consequences is the existence of hyperbolic
periodic orbits in $K$ of all periods $kn, n\ge 1$.  Note, however,
that $k$ is not specified in the horseshoe theorem, and in most cases
there is no way to estimate it (but see Mischaikow \& Mrozek
\cite{MM95}).  Collins \cite{Co99} has shown that a
topologically transverse homoclinic point implies the existence of
periodic points of all sufficiently high minimum periods; estimating
such periods, however, requires detailed knowledge of the associated
homoclinic tangle.

While the horseshoe theorem guarantees infinitely many periodic
orbits, it is insufficient for the existence of a second fixed
point. For example, the toral diffeomorphism induced by the matrix $
\left[\begin{smallmatrix} 2 & 1\\ 1 & 1\end{smallmatrix}\right]$ has
only one fixed point, even though transverse homoclinic points are
dense.

It turns out that for diffeomorphisms of the plane, even a
nontransverse homoclinic point implies a second fixed point; in fact
there is a block of fixed points having index $+1$.  But the proof of
this (Hirsch \cite{Hi99}), based on Brouwer's Plane Translation
Theorem, gives no indication of the location of such a block.

In this paper we consider a saddle fixed point $p$ for an orientation
preserving homeomorphism $f$ of a surface $X$.  Let $p'$ be a
homoclinic point associated to $p$, i.e., a point different from $p$
wherethe stable and unstable curves $W_s (p)$, $W_u (p)$ meet; no
transversality or even crossing of these curves is assumed.  Suppose
$\Gamma=J_s\cup J_u$ is a homoclinic loop $p$, where $J_s$ and $J_u$ are
arcs in $W_s (p)$ and $W_u (p)$ respectively, having common endpoints
$p,p'$.  Assume there is a closed 2-cell, with interior $V$, whose
boundary is the union of two arcs in $W_s (p)$ and $W_u (p)$
respectively, having endpoints $p$ and $p'$ in common but otherwise
disjoint; such a 2-cell always exists when $X$ is simply connected.

Our main result, stated more precisely in Theorem \ref{th:A}, is this:
 
\smallskip
\noindent{\em If $\Lambda$ is a Jordan curve bounding an open 2-cell
$V$, there exists $\rho\in\{1,2\}$ such that the fixed point index of
$f^n$ in $V$ for all $n\ne 0$, and $\rho$ depends only on the geometry
of $V$.}

\medskip
\noindent
An immediate consequence is that for every $n\ge 2$, every map
sufficiently close to $f^n$ has a block of fixed points
in $V$ of index $\rho$ .

Theorem \ref{th:E} is a similar result for homoclinic loops that are
homotopically trivial, but not necessarily Jordan curves.

The main theorem is stated in 
Section \ref{sec:statement}, and  several applications are derived.
Section \ref{sec:f} contains the proof of the main theorem. 
\section{The main result and applications}   \mylabel{sec:statement}
$\ZZ, \NN $ and $\Np$ denote the integers, natural numbers and
positive natural numbers.

All maps are assumed continuous; $\approx$ denotes homeomorphism. 

For any map $g$, the maps $g^n, n\ge 1$ are defined recursively by
$g^1=f$ and $g^{n+1} (x)=g (g^n (x))$ provided $g^n(x)$ is in the
domain of $g$.

$X$ always denotes a connected, oriented surface with metric $d$, and
$f\co X\to S$ is an orientation preserving injective map.  We call $f$
a diffeomorphism when $f$ and $f^{-1}$ are $C^1$ (continuously
differentiable). 

The {\em orbit} of $x$ is the set $\ga x =\{f^i (x):i\in \ZZ\}$.  The
fixed point set of $f$ is denoted by $\Fix f$.  We call $q\in \Fix
(f)$  {\em smooth} if it belongs to a coordinate chart
in which $f$ is represented by a $C^1$ map; such a chart is 
{\em smooth} for $p$.  If $f$ is $C^1$, of
course all fixed points are smooth.  But in many constructions 
some fixed points of a  nonsmooth maps are smooth, as when a
diffeomorphism of the plane is extended to the 2-sphere.  

Let $q\in\Fix f$ be smooth.  We call $q$ 
{\em simple} if $1$ is not an eigenvalue of the linear operator
$df_q$, {\em hyperbolic} if no eigenvalue lies on the unit circle
$S^1\subset \bf C$, a {\em sink} if eigenvalues are inside $S^1$, a
{\em source} if they are outside, and {\em elliptic} if the
eigenvalues are on $S^1$ but different from $1$. 

A fixed point $p$ is a {\em saddle} if it is not in the boundary of
$S$, and there is a chart at $p$ in which $f$ is locally
represented as a linear map $\left[\begin{smallmatrix} \mu & 0\\ 0 &
\lam \end{smallmatrix}\right]$, and either $\mu>1>\lam>0$, making $p$
a {\em direct} saddle, or $\mu<-1<\lam<0$, defining a {\em twisted
saddle}.  Such a chart is {\em diagonalizing}.  By the
Hartman-Grobman linearization theorem (Hartman \cite{Ha64}), for $p$
to be a saddle it is sufficient for there to be a smooth chart at $p$
in which $df_p$ has eigenvalues $\mu, \lam$ as above.

An {\em $n$-periodic} point for $f$ means a fixed point $z$ for $f^n,
\:n\ge 1$.  When $n$ is the minimum period, $\ga z$ is an {\em
$n$-orbit}.  An $n$-periodic point is simple, hyperbolic, and so
forth, when it has the corresponding property as a fixed point for
$f^n$.

The {\em stable curve} $W_s=W_s (p)$ of a saddle fixed point $p$ is
the connected component of $p$ in the set of $x$ for which there is a
convergent sequence $x_k\to p$ in $X$ with $x_0=x$ and
$f(x_k)=x_{k+1}$.  The {\em unstable curve} $W_u$ at $p$ is defined as
the stable curve for $f^{-1}$.  Note that $W_s$ and $W_u$ are mapped
homeomorphically onto themselves by $f$.
Owing to the linearization assumption, there are bijective maps
$\zeta_u, \zeta_s:\RR \to W_s$ taking $0$ to $p$, called {\em
parametrizations} of $W_u, W_s$ respectively.  The images of
$[0,\infty)$ and $(-\infty,0]$ are the four {\em branches} at $p$.

A {\em homoclinic point} for $p$ is any point $p'\in W_s\cap
W_u\setminus\{p\}$, in which case the {\em homoclinic loop} $\Lam$
defined by $p'$ is the closed path formed by the two arcs $J_s \subset
W_s,\:J_u\subset W_u$ having common endpoints $p$ and $p'$.  There
corresponds an element $[\Lam]$ of the fundamental group of $X$ at
$p$, determined by first traversing $\Lam$ from $p$ to $p'$ in $J_u$
and then from $p'$ to $p$ in $J_s$.  If $[\Lam]$ is the unit element,
$\Lam$ is a {\em inessential homoclinic loop}.  $\Lambda$ is {\em
simple} if $J_u\cap J_s =\{p,p'\}$, in which case $\Lambda$ is
homeomorphic to the unit circle.  Every homoclinic
loop contains a simple homoclinic loop.

Suppose $\Lam$ is a simple homoclinic loop in $X$ bounding a closed
2-cell $D\subset X$.  The corresponding open 2-cell $V=D\setminus
\partial D$ is a {\em homoclinic cell}.  We call $V$ is a {\em
positive} cell provided some diagonalizing chart takes $p$ to the
origin $0\in\R 2$ and a neighborhood of $p$ in $D$ onto a neighborhood
of $0$ in the first quadrant.  In the contrary case $D$ is a {\em
negative} cell: there is a diagonalizing chart taking a neighborhood
of $p$ in $D$ onto a neighborhood of the origin in the complement of
the open first quadrant (see Figure 1).  Thus when seen through a
diagonalizing chart, a positive cell appears convex near $p$, while a
negative region appears concave.
\begin{figure} [ht] 
\centering
\epsfig{figure=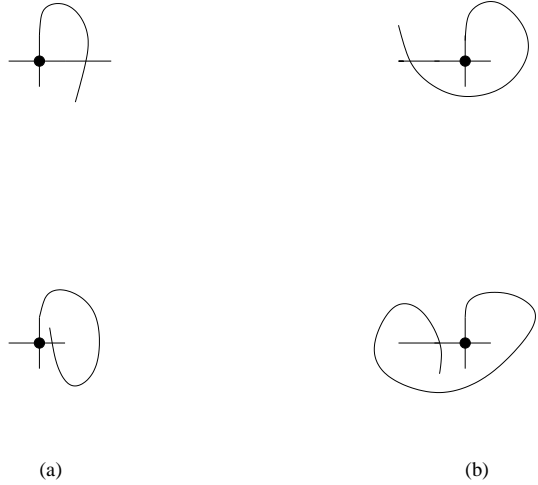, height=2.5in} 
\caption
{
\small Homoclinic cells:\hspace{.30cm} (a) positive \hspace{.25cm} (b)
negative 
\label{fig1}
}
\end{figure}

Let $U\subset X$ be an open set such that $U\cap\Fix f$ is compact.
The {\em fixed point index of $f$ in $U$} is denoted by $I(f,U)\in
\ZZ$; if it is nonzero, there exists a fixed point in $U$ (Dold
\cite{Do65}).
When $U$ is a coordinate chart identified with an open set in $\R 2$,
we can calculate $I(f,U)$ as follows.  Let $M\subset U$ be a compact
surface with boundary whose interior contains $\Fix f \cap U$.  Then
$I(f,U)$ is the degree of the map
\[
  \partial M \to S^1,\quad x\mapsto \frac {x-f(x)}{\vert\vert
  x-f(x)\vert\vert},
\] 
where $\partial M$ and $S^1$ inherit their orientations from $\R 2$.
If $\partial M$ is replaced by any oriented Jordan curve $\Gamma$ on
which $f$ has no fixed points, the same formula defines the {\em index
of $f$ along $\Gamma$}.

Let $B\subset X$ be a {\em block} of fixed points, i.e., $B$ is
compact and relatively open in $\Fix f$.  There exists an open
neighborhood $U_0\subset X$ such that $B=\Fix f \cap U_0$.  The number
\[\Ind (f,B)=\Ind (f|U_0)\in\ZZ,\] 
called the {\em index of $f$ at $B$}, is independent of the choice of
$U_0$. When $p$ is an isolated fixed point we set $\Ind
f,\{p\})=I(f,p)$, called the {\em index of $f$ at $p$}.  A direct
saddle has index $-1$. Twisted saddles, sources, sinks and elliptic
fixed points have index $+1$.

The following assumptions are in force throughout the rest of
this article:
\begin{hypothesis} \mylabel{th:standing}
{~}
\begin{itemize}

\item $f\co X\approx X$ is an orientation-preserving homeomorphism
of a surface $X$.

\item
$p\in X\setminus \partial X$ is a direct saddle fixed point for $f$.

\item  $V\subset X$ is an open 2-cell bounded by a simple homoclinic loop
$\Lambda$ at $p$.
\end{itemize}
\end{hypothesis}
To $V$ we assign the number
\[\rho =\rho(V)=\left\{ \begin{array}{ll} 1 & \mbox{ if $V$ is a positive
		region } \\ 2 & \mbox{ if $V$ is a negative region.}
		\end {array} \right\}
\]
For  each $n\in \Np$ we define an open set $V_n\subset V$,
\[V_n =V_n (f)=\{x\in V:f^i (x)\in V,\: i=1,\ldots,n-1\}\]
Thus  $\Fix {f^n}\cap V_n$ is the union of the $n$-periodic
orbits in $V$.  

\smallskip
The following is our fundamental result:
\begin{theorem}		\mylabel{th:A}
$\Fix {f^n}\cap V_n$ is a block of fixed points for $f^n$ of index
$\rho (V)$, for all $n\ge 1$.
\end{theorem}

\noindent
Before giving the proof of Theorem \ref{th:A} in Section \ref{sec:f},
we present several consequences. Hypothesis \ref{th:standing} is
always assumed.

\subsection{Homeomorphisms of the sphere}
Assume $g\co S^2\approx S^2$ is an orientation preserving
homeomorphism having a simple homoclinic loop $\Lam$ at a direct
saddle.
\begin{theorem}		\mylabel{th:6}
The fixed point index of $g$ in one
of the two complementary  components 
of $ \Lam$ is $1$, and the index in the other  is $2$.  
\end{theorem}
\begin{proof}
This follows from Theorem \ref{th:A}, because one complementary
component of $\Lam$ has positive type and the other has negative type.
\end{proof}

The persistence of blocks having nonzero index implies:
\begin{corollary}		\mylabel{th:8}
Every map $S^2\to S^2$ sufficiently close to $g$ has at least 3
fixed points.
\end{corollary}

A homoclinic loop constrains fixed point indices.  Suppose for example
that there are exactly 3 fixed points: a direct saddle and two other
fixed points with respective indices 5 and -2.  Then the
saddle does not admit a homoclinic point.

\subsection{Inessential homoclinic loops and Nielsen classes}

Fixed points $a, b$ are in the same {\em Nielsen class} provided they
are endpoints of a path that is homotopic to its composition with $f$,
keeping endpoints fixed.  Equivalently, $f$ is covered by a map in a
universal covering space having fixed points over $a$ and $b$.

When $X$ is compact, every Nielsen class is a block of fixed points,
and the {\em Nielsen number} of $f$ is the number of Nielsen classes
having nonzero index.  This number, a homotopy invariant of $f$, is a
lower bound for the number of fixed points for any map homotopic to
$f$.

\begin{theorem}		\mylabel{th:E}
Assume $p$ belongs to an inessential homoclinic loop.  Then its Nielsen
class contains a block of positive index, and such a block must
contain a fixed point $q\ne p$.  When the Nielsen class of $p$ is
finite, $q$ can be chosen with positive index.
\end{theorem}
Thus in the presence of a inessential homoclinic loop, the number of
fixed points exceeds the Nielsen number.  Theorem \ref{th:C} is a
similar result for Lefschetz numbers.

\begin{corollary}		\mylabel{th:Ecor}
If a direct saddle $p$ is the only member of its Nielsen class, then
$p$ does not belong to an inessential homoclinic loop.
\end{corollary}
  
In proving Theorem \ref{th:E} we assume $X$ is not simply connected,
otherwise using Theorem \ref{th:A}.  As $X$ is orientable, there is a
universal covering space $\pi:\R 2\to X$.
Choose $\til p\in\pi^{-1}(p)$, and let $\til f\co \R 2\to\R 2$
be the unique lift of $f$ with a fixed point at $\til p$.  Then $\til
p$ is a direct saddle for $\til f$.

Let $\Gam$ be a null homotopic homoclinic loop at $p$.  There is a
unique homoclinic loop $\til \Gam \subset \R 2$ for $\til f$ that
contains $\til p$ and projects onto $\Gam$ under $\pi$.  Let
$\til\Lam\subset \til \Gam$ be a simple homoclinic loop at $\til p$.
There is a unique an open 2-cell $V\subset \R 2$ bounded by $\Lam$,
and $\ov V$ is a closed 2-cell.  Applying Theorem \ref{th:A}, we choose a
block $L\subset \Fix {\til f}\cap V$ such that
\[\Ind(\til f, L)=\sigma\in \{1,2\}\]
Notice that $\pi (L)$ lies in the Nielsen class of $p$. 

Every fixed point $\til z\in \pi^{-1}(p)$ has index $-1$, since $\til
f$ in a neighborhood of $\til z$ is conjugate to $f$ in a neighborhood
of $\pi (z)$.  Because $\ov V$ is compact, $\pi^{-1}(p)\cap V$ is
finite.  Therefore $L\setminus \pi^{-1} (p)$ is nonempty, for
otherwise $L$ would be a nonempty finite subset of $\pi^{-1} (p)$ and
thus have negative index.

It follows that $\pi(L\setminus \pi^{-1}(p))$ is a nonempty subset of
$\Fix f$ disjoint from $p$, contained in the Nielsen class of $p$.
Suppose this class is finite.  Then $L\setminus \pi^{-1}(p))$ is
finite.  Let
$ L\cap \pi^{-1}(p) $ have cardinality $\nu,\;1\le \nu<\infty$.
Then
\[
\begin{split}
\Ind (\til f, L\setminus \pi^{-1}(p))) &=\Ind (\til f, L)-
\Ind (\til f, L\cap (\pi^{-1}(p)))\\
&=\Ind (\til f, L)-\nu \,\Ind (f,p))\\
&= \sig +\nu \quad\ge\quad 2
\end{split}
\]
Therefore exists $\til q\in L\setminus \pi^{-1}(p)$ with $0<\Ind (\til
f, \til q) =\Ind (f, \pi (\til q))$, and $\pi (\til q)$ is in the
Nielsen class of $p$.  This completes the proof of Theorem \ref{th:E}.

\subsection{Periodic orbits in a homoclinic cell}
The following theorem can be used to demonstrate the existence of
infinitely many periodic orbits in situtations where the horseshoe
theorem may not apply:
\begin{theorem}		\mylabel{th:A1}
Let $r\in\NN$ be such that every
$2^k$-orbit, $0\le k\le r$ in the homoclinic cell $V$ is hyperbolic.
Then  either: 
\begin{description}
\item [(a)] 
$V$ contains an
attracting or repelling $2^k$-orbit for some $k\in\{0,\dots,r\}$,
\end{description}
 or else
\begin{description}
\item [(b)]  $V$ contains a twisted saddle orbit of
cardinality $2^k$ for every $k=0,\dots,r$. 
\end{description}
\end{theorem}
\begin{proof}
Suppose (a) does not hold.  Fix $n=2^k, \:0\le k\le r$ and let
$B\subset \Fix {f|V_n}$ be a block having index $\rho\in\{1,2\}$
(Theorem \ref{th:A}).  Then some $q\in B$ has index $1$ for $f^n$.
Since (a) is ruled out, $q$ is not a source or sink for $f^n$.  The
only other possibility for a hyperbolic, index $1$ fixed point for
$f^n$ is a twisted saddle.  This implies $n$ is the minimal period for
$q$.  Thus (b) holds.
\end{proof}

\begin{corollary}		\mylabel{th:A2}
Assume every periodic orbit in $V$ whose cardinality is a power of $2$
is a saddle.  Then $V$ contains a twisted saddle orbit of cardinality
$2^k$ for every $k\in\NN$.
\end{corollary}

\begin{corollary}		\mylabel{th:A3}
If $f$ is $C^1$ and $0<\Det df_x <1$ in a dense subset of $V$ and all
periodic points in $V$ are hyperbolic, then $V$ contains either a
periodic attractor, or an orbit of cardinality $2^k$ for every
$k\in\NN$
\end{corollary}

It is interesting to compare these results to a theorem of Franks
\cite{Fr77}.  Specialized to an orientation-preserving diffeomorphism
of the 2-sphere, it states:
\begin{quote}
{\em If all periodic points are hyperbolic, and at
most one orbit whose cardinality is a power of 2 is repelling or
attracting, then there are infinitely many periodic orbits.}
\end{quote}
Corollary \ref{th:A2} makes no assumptions on orbits outside the
homoclinic cell $V$, but does not allow any attractors or repellors of
cardinality $2^k$ in $V$.  It gives sharper information than the
conclusion Franks' theorem on the periods and locations of periodic
orbits.

It is not trivial to construct diffeomorphisms of the disk or sphere,
all of whose periodic orbits are saddles; but examples are known that
even have the Kupka-Smale property, i.e., stable and unstable curves
of periodic points have only transverse intersections (Bowen \& Franks
\cite{BF76}, Franks \& Young \cite{FY80}).  Gambaudo {\em et.\ al}
\cite{GVT89} construct real analytic Kupka-Smale examples on the disk.

\subsection{Lefschetz numbers}
Let $\#Q$ the cardinality of a set $Q$.

Suppose $X$ is a compact surface and $h:X\to X$ is continuous.  The
{\em Lefschetz number} $\mathsf {Lef} (h)$ is the alternating
sum of the traces of the induced endomorphisms of the singular
homology groups $H_i (X), \:i=0,1,2$; it  equals $\Ind (h, X)$.  Lefschetz
proved that when the fixed point set is finite, $\mathsf {Lef} (h)$ is
the sum of the fixed point indices.  When every fixed point has index
$+1$, $-1$ or $0$ this gives the useful estimate
\[\#\,\Fix h \ge |\mathsf {Lef} (h)|\]
The following results show that when fixed points are simple,
homoclinic cells entail the existence of more fixed points than are
counted by the Lefschetz number.
\begin{theorem}		\mylabel{th:C}
Assume $X$ is a compact surface, $\Fix f$ is finite, and every fixed
point has index $+1$, $-1$ or $0$.  If $f$ admits a
homoclinic cell, then
\[\#\,\Fix f\:\ge\: |\mathsf {Lef}(f)\,+1-\rho| +1+\rho \:
\ge \: |\mathsf {Lef}(f)|+2
\] 
\end{theorem}
\begin{proof}
For any open set $A\subset X$, summing indices over fixed points $z\in
A$ gives:
\[
\begin{split}
|\Ind (f,A)| &=\left| \sum_{z} \Ind (f,z)\right|
\le \sum_{z} |\Ind (f,z)|\\[1ex]
& \le \#\,(\Fix f\cap A)
\end{split}
\]  
Applying this to a homoclinic cell $V$, from Theorem \ref{th:A} we get
\[
\begin{split}
\#\,(\Fix f\cap\ov V) &=1+\#\,(\Fix f\cap V) \ge 1+|\Ind (f,V)|\\
&= 1+\rho
\end{split}
\]
because $\Ind (f, V)=\rho$. 
Also
\[
\begin{split}
\#\,(\Fix f\cap (X\setminus \ov V)) 
	&\ge |\Ind (f,X\setminus \ov V )|\\
	&=|\Ind (f, X)- (\Ind (f,p)+\Ind (f, V))|\\
	&=|\mathsf {Lef}(f) +1 -\rho |
\end{split}
\]
because $\Ind (f,p)=-1$. 
Therefore
\[
\begin{split}
\#\,\Fix f 
&=\#\,(\Fix f\cap  \ov V)\,+\,\#\,(\Fix f\cap (X\setminus \ov V))\\ 
&\ge (1+\rho)\, +\, |\mathsf {Lef}(f) +1 -\rho|\\
&\ge   |\mathsf {Lef}(f)| + 2
\end{split}
\]
\end{proof}

\begin{corollary}		\mylabel{th:Ccor}
Assume $X$ is a compact surface, $\Fix f$ is finite, and every fixed
point has index $+1$, $-1$ or $0$.  If $\#\,\Fix f\le |\mathsf
{Lef}(f)| +1$, there are no homoclinic cells.
\end{corollary}

\section{Fixed point indices and retractions} 
\mylabel{sec:f}
This section contains the proofs of Theorems \ref{th:A} and
\ref{th:E}.  Hypothesis \ref{th:standing} continues to hold.  Let
$D\subset X$ denote the closure of the
homoclinic cell $V$.  Then $D$ is a compact 2-cell whose boundary is
the simple homoclinic loop $\Lambda$. 

A {\em retraction} of a space $Y$ onto a subset $Y_0\subset Y$ is a
map $Y\to Y_0$ fixing every point in $Y_0$.

\begin{lemma}		\mylabel{th:0}
Assume we are given $n\in\Np$ and a map $g\co D\to D$ with the
following properties:
\begin{description}

\item [(i)] $g$ coincides with $f^n$ on a neighborhood of $p$ in $D$

\item [(ii)]  $\Fix {g}=K\cup\{p\}$
where $K\subset V$ is compact.
\end{description}
Then  $\Ind (g, K)=\rho (V)$
\end{lemma}
\begin{proof}
Fix a coordinate chart in which $p$ is the origin and $f^n$ is
represented by a linear map
\[T (x,y)= (\lam x, \mu y),\quad 0<\lam<1<\mu\] 
We identify points near $p$ with their images in $\R 2$ under this
chart.

Consider the case that $V$ is a positive homoclinic cell ($\rho=1$).
Then there is a compact disk neighborhood $N\subset \R 2$ centered at
the origin, meeting $D$ only in one of the four closed quadrants; to
fix ideas, we assume it is the first quadrant $Q_I$.  We take $N$ so
small that $g$ coincides with $T$ in $N\cap D\:$, $N\cap
K=\varnothing\:$, and $N\cup D$ is a 2-cell.

Choose a retraction $s:N\to N\cap Q_I$.  We compute the fixed point
index $\Ind (T\circ s, N)$.  Let $\eps>0$ be so small that the disk
$D_\eps$ of radius $\eps$ lies in $N$.  Let $S^1_\eps$ denote the
circle bounding $D_\eps$.  Since $T\circ s$ has the unique fixed point
$0$, the index equals the degree of the map
\[u:S^1_\eps \to S^1,\: 
z\mapsto
\frac{z-T\circ s(z)}
 {|| z-T\circ s
(z)||}\]
The retraction $s$ sends any point $z\in N\setminus Q_I$ to the unique
point $s(z)\in \partial Q_I$ such that $z$ and $s(z)$ are the
endpoints of a line segment having slope $1$; and $s$ is the identity
on $N\cap Q_I$.  A simple computation shows that $u$ takes no values
in the first quadrant of the unit circle, and thus has degree zero.
Thus $\Ind (T\circ s, N)=0$.

Now consider the map $h:N\cup D\to D \subset N\cup D$ defined to be
$T\circ s$ in $N$ and $g$ in $D$; this definition is consistent
because $s$ is a retraction and $g$ coincides with $T$ in $N\cap D$.
Clearly
\[\Fix {h}= \{p\}\cup K \subset \Int (N\cup D)\]
Therefore 
\[
  \mathsf{Lef}(h)=\Ind (h,\Int (N \cup D))=\Ind (h,p)+\Ind (h, K)
\]
Note that $\mathsf{Lef}(h)=1$ because $N \cup D$ is a compact 2-cell,
and $\Ind (h, p)=\Ind (T\circ s, N)=0$.  Hence
\[1=\Ind (h, K)=\Ind (g, K)\]
as required. 

When $V$ is a negative homoclinic cell, we can assume $N\cap D$ {\em
excludes} the interior of the first quadrant.  The retraction $s:N\to
N\setminus \Int Q_I$ is defined by sending $z\in N\cap Q_I$ to the
unique point of $\partial Q_I$ such that $z$ and $r(z)$ are the
endpoints of line segment having slope $1$; and $r$ is the identity on
$N\setminus Q_I$.  The degree of $T\circ r$ in this case is $-1$.  Define $h$
as above.  An argument similar to the preceding shows that
\[
  \begin{split}
  1 &=\Ind (h, {\rm Int}(N\cup D))=\Ind (h,\{p\})+\Ind (h,{\rm
       Int}(N\cup D) )\\
  &=-1+\Ind (g,{\rm Int}(N\cup D)) 
\end{split}
\]
\end{proof}

Let $J_u \subset W_u(p)$ and $J_s \subset W_s(p)$ denote the two
compact arcs whose union is $\Lambda$; these arcs meet at their common
endpoints, which are $p$ and the homoclinic point $p'\ne p$, but
nowhere else.

Our next goal is the following result:

\begin{proposition}		\mylabel{th:r}
There is a retraction 
\[r:f(D)\cup D\to D\]
such that 
\begin{equation}		\label{eq:r}
r(f(D)\setminus D)\subset J_s
\end{equation}
\end{proposition}
\begin{proof}
We first prove
\begin{equation}		\label{eq:a}
J_u \cap \clos ( f (D)\setminus D )=\{p,p'\}
\end{equation}
or equivalently,
\[J_u \cap \clos(f (V)\setminus D)=\{p,p'\}\]
Suppose (\ref{eq:a}) is false, so that there exists 
\[b\in J_u\setminus\{p,p'\} \cap \clos(f (V)\setminus D)\]
Then $b=\lim_{i\to\infty} f(a_i)$ for some sequence $a_i\in V\setminus
f^{-1}D$, and $b=f(a)$ by continuity.  I claim $f$ maps a relatively
open neighborhood $N_a\subset D$ of $a$ onto a relatively open
neighborhood $f (N_a))\subset f (D)$ of $f(a)$.  This is because $f$
maps the interior of $D$ onto the interior of $f (D)$.  The assumption
that $p$ is a direct saddle implies $f$ preserves orientation, and
$f^{-1}|J_u$ preserves orientation in $J_u$.  From this it follows
that $N_a$ and $f(N_a)$ abut $J_u$ from the same side.  Consequently
$f(N_a)$ contains a relatively open neighborhood $N_b\subset D$ of
$b$.  For sufficiently large $i$ we have $a_i\in f^{-1}N_b$ and thus
$a_i\in f^{-1}D$.  This contradiction completes the proof of
(\ref{eq:a}).

From Equation (\ref{eq:a}) we see that
\begin{equation}		\label{eq:b}
\clos (f(D)\setminus D)\cap\ D \subset J_s
\end{equation}
Note also
\[
\begin{split}
f(D)\cup D	 &=  \clos (f(D)\setminus D) \cup D,\\
\clos (f(D)\setminus D) \cap D =& \clos (f(D)\setminus D)\cap\partial D
\subset J_s
\end{split}
\]
By Tietze's extension theorem there is a retraction 
\[r_0:\clos (f(D)\setminus D)\cup J_s \to J_s\]
and $r_0$ agrees with the identity map of $D$ on the intersection of
their domains, which by (\ref{eq:b}) is $J_s$.  Thus $r_0$ and the
identity map of $D$ fit together to give the desired retraction $r$.
\end{proof}

From now on $r:f(D)\cup D\to D$ denotes a retraction as in Proposition
\ref{th:r}.

\begin{lemma}		\mylabel{th:3a}
Let $n\in \NN$.  For every $q\in\Fix {f^n}\cap V_n$  there is a
neighborhood $U\subset V_n$ of $q$ such that $f^n|U= (r\circ f)^n|U$.
\end{lemma}
\begin{proof}
The definition of $V_n$ implies $f^j (q)\in V_n\subset V$ for all
$j\in \NN$.  Therefore $q$ has a neighborhood $U$ such that $f^i
(U)\subset V_n$ for $i=0,\dots, n$.  Assume inductively that $0\le
i<n$ and $f^i|U= (r\circ f)^i|U$; the case $i=0$ is trivial.  For
$x\in U$ we have $(r\circ f)^i (x)=f^i(x)$, and both $f^i(x)$ and
$f^{i+1} (x)$ are in $V$ because $x\in V_n$.  Hence
\[(r\circ f)^{i+1} (x)=(r\circ f) (f^i (x))=r (f^{i+1} (x))=f^{i+1}
(x)\]
because $r$ and $f$ coincide on $V$.  This completes the induction. 
\end{proof}
\begin{lemma}		\mylabel{th:3b}
$\Fix {(r\circ f|D)^n}=\{p\}\cup(\Fix {f^n}\cap V_n)$ for all $n\ge 1$.
\end{lemma}
\begin{proof} 
Let $x\in D\setminus\{p\}$ be $n$-periodic for $r\circ f$.  We first
show $x\notin J_s$.  We know that $J_s$ is invariant under $f$, and
$r|J_s$ is the identity because $J_s\subset D$.  Thus $r\circ f|J_s$
coincides with $f|J_s$, whose only periodic point is $p$.  The
foregoing implies no point on the orbit $x$ under $r\circ f$ lies in
$J_s$.  Therefore no point $y$ in this orbit maps outside $D$ under
$f$, for otherwise $(r\circ f) (y)\in J_s$ by Equation (\ref{eq:r}).
This proves $\ga x \subset D$, and an induction that $(r\circ f)^k x=
f^k x$ for all $k$.  Since $J_u\setminus p$ is contains no periodic
points for $f$, the conclusion follows.
\end{proof}

\subsection{Proof of Theorem \ref{th:A}} 
The set $B=\Fix {f^n}\cap V_n$  is open in $\Fix {f^n}$ because $V_n$
is open. 
We prove $B$  compact by showing it is closed in $D$.  Since
$\ov B\cap\partial D\subset \{p\}$,  it suffices to prove that $p$
is not a limit point of $B$.  Clearly $p\notin B$, and $p$, being a
saddle, has a neighborhood in which the only point of period $n$ is
$p$.  Therefore $B$ is a block.  

To prove $\Ind (f^n, B)=\rho$, let $r:f(D)\cup D\to D$ be a retraction
as in Proposition \ref{th:r}.  Lemmas \ref{th:3b} and \ref{th:3a} show
that $\Ind (f^n, B)=\Ind ( (r\circ f|V)^n, B)$.  Now apply Lemma
\ref{th:0} to $g:=(r\circ f|D)^n$ to conclude that $\Ind ((r\circ
f|V)^n, B)=\rho$. \qed

\end{document}

%% file: preamble4.tex
\def\bc{\begin{center}}
\def\ec{\end{center}}

\newtheorem{theorem}{Theorem}[section]
\newtheorem{lemma}[theorem]{Lemma}
\newtheorem{corollary}[theorem]{Corollary}
\newtheorem{proposition}[theorem]{Proposition}

{\theorembodyfont{\rmfamily}

}
{\theoremstyle{break}	
{\theorembodyfont{\rmfamily} 
				
}}
{\theoremstyle{break}	
			\newtheorem{hypothesis}[theorem]{Hypothesis}
}

{\end{enumerate}}

\newenvironment{proof}[1]{\smallskip \noindent {\bf #1}}{\qed\smallskip}

\def\qed{\ifhmode\unskip\nobreak\fi\ifmmode\ifinner\else\hskip5pt\fi\fi
 \hfill\hbox{\hskip5pt\vrule width4pt height6pt depth1.5pt\hskip1pt}}

\newcommand{\ov}[1]{\ensuremath{\overline{#1}}} 
 
\newcommand{\Ind}{\ensuremath{\operatorname{\mathsf {Ind}}}}
\newcommand{\Int}{\ensuremath{\operatorname{Int}}}
\newcommand{\Det}{\ensuremath{\operatorname{Det}}}
\newcommand{\clos}{\ensuremath{\operatorname{clos}}}

\newcommand{\RR}{\ensuremath{{\bf R}}}     
\newcommand{\R}[1]{\ensuremath{{\bf R}^{#1}}} 
\newcommand{\ZZ}{\ensuremath{{\bf Z}}}	   
\newcommand{\NN}{\ensuremath{{\bf N}}}	   
\newcommand{\Np}{\ensuremath{{\bf N}_{+}}}  

\def\d#1dt{\frac{d#1}{dt}}    
\newcommand{\eps}{\ensuremath{\epsilon}}
\newcommand{\lam}{\ensuremath{\lambda}}
\newcommand{\Lam}{\ensuremath{\Lambda}}

\newcommand{\sig}{\ensuremath{\sigma}}

\newcommand{\Gam}{\ensuremath{\Gamma}}

\newcommand{\Fix}[1]{\ensuremath{\mathsf {Fix}\,(#1)}}

\def\ga#1{{\gamma}(#1)}
\def\Fix#1{\mathsf {Fix}\,(#1)}
\def\til{\tilde}

\newcommand{\co}{\colon\thinspace} 
 